\documentclass{article}
\usepackage{amssymb}

\usepackage{graphicx}
\usepackage{amsmath}

\newtheorem{theorem}{Theorem}

\newtheorem{conjecture}[theorem]{Conjecture}

\newtheorem{lemma}[theorem]{Lemma}

\newtheorem{remark}[theorem]{Remark}

\newenvironment{proof}[1][Proof]{\textbf{#1.} }{\ \rule{0.5em}{0.5em}}
\input{tcilatex}

\begin{document}

\title{\bf{Parabolic twists for algebras $sl(n)$.}}
\author{Vladimir Lyakhovsky \\
Department of Theoretical Physics,\\
Sankt-Petersburg State University,\\
Ulianovskaya 1, Petrodvoretz, 198904 \\
Sankt Petersburg, Russia}
\maketitle

\begin{abstract}
New solutions of twist equations for universal enveloping algebras
$U\left( A_{n-1}\right) $ are found. They can be presented as products of
full chains of extended Jordanian twists $\mathcal{F}_{\widehat{ch}}$, Abelian
factors ("rotations") $\mathcal{F}^{R}$ and sets of quasi-Jordanian twists
$\mathcal{F}^{\widehat{J}}$. The latter are the generalizations of Jordanian
twists (with carrier $b^2$) for special deformed extensions of the Hopf algebra
$U\left( b^2 \right)$. The carrier
subalgebra $g_{\mathcal{P}}$ for the composition $\mathcal{F}_{\mathcal{P}}=%
\mathcal{F}^{\widehat{J}}\mathcal{F}^{R}\mathcal{F}_{\widehat{ch}}$ is a
nonminimal parabolic subalgebra in $A_{n-1}$, $g_{\mathcal{P}%
}\cap \mathbf{N}_{g}^{-}=\emptyset $. The parabolic twisting elements
$\mathcal{F}_{\mathcal{P}}$ are obtained in the explicit form. The details
of the construction are illustrated by considering the examples $n=4$ and
$n=11$.
\end{abstract}

\section{Introduction}

The parabolic subgroups play an essential role in group theory and
applications. Parabolic subgroups of a real simple Lie groups can be treated
as the groups of motion of noncompact symmetric spaces with the special
properties of geodesics \cite{ONI}. The quantized versions of such subgroups
are important for the models of noncommutative space-time \cite{CKNT,WESS}.
Quantizations of antisymmetric $r$-matrices, (solutions of classical
Yang-Baxter equation) form an important class of deformations producing the
quantized groups. The quantum duality principle \cite{Drin-3,Sem-Tyan}
provides the possibility to describe quantum groups with dual Lie algebras
in terms of quantized universal enveloping algebras. Thus to construct a
quantized parabolic subgroup of a simple Lie group it is sufficient to
consider a parabolic subalgebra $g_{\mathcal{P}}$ in a simple Lie algebra $g$
and to develope for the universal enveloping algebra $U\left( g\right) $ a
quantization where the parabolic properties are selfdual. The important
class of such quantum deformations is known as twisting \cite{DRIN-83}.

If we supply the elements of Lie algebra $g$ with primitive coproducts $%
\Delta ^{\left( 0\right) }$ and consider $U(g)$ as a Hopf algebra with the
costructure generated by $\Delta ^{\left( 0\right) }$ then the invertible
solution $\mathcal{F}\in U(g)^{\otimes 2}$ of the twist equations \cite
{DRIN-83}
\begin{equation}
\begin{array}{l}
\left( \mathcal{F}\right) _{12}\left( \Delta ^{\left( 0\right) }\otimes
\mathrm{id}\right) \mathcal{F}=\left( \mathcal{F}\right) _{23}\left( \mathrm{%
id}\otimes \Delta ^{\left( 0\right) }\right) \mathcal{F}, \\
\left( \epsilon \otimes \mathrm{id}\right) \mathcal{F}=\left( \mathrm{id}%
\otimes \epsilon \right) \mathcal{F}=1.
\end{array}
\label{twist-equations}
\end{equation}
allows to transform $U(g)$ into the Hopf algebra $U_{\mathcal{F}}(g)$,
\begin{equation}
\mathcal{F}:U(g)\longrightarrow U_{\mathcal{F}}(g).
\end{equation}
The twisted Hopf algebra $U_{\mathcal{F}}(g)$ has the same multiplication,
unit and counit as $U(g)$ but the twisted coproduct and antipode given by:
\begin{equation}
\Delta _{\mathcal{F}}(a)=\mathcal{F}\Delta ^{\left( 0\right) }(a)\mathcal{F}%
^{-1},\quad S_{\mathcal{F}}(a)=vS(a)v^{-1},
\end{equation}
\begin{equation*}
v=\sum f_{i}^{(1)}S(f_{i}^{(2)}),\qquad a\in U(g).
\end{equation*}
The minimal subalgebra $g_{c}\subset g$ necessary for the definition of the
twist $\mathcal{F}$ is called the \textbf{carrier} of the twist \cite{GERST}.

To investigate the deformation quantization problems and in particular
to apply the twisted parabolic groups to this study
explicit solutions of twist equations are necessary \cite{BGGS}.
Only certain classes of twists $\mathcal{F}$ are known in the explicit
form \cite{RESH, OGI, KLM,KLO, AKL}. All such solutions are factorizable,
i. e. can be presented as a product of the so
called \textbf{basic twisting factors} \cite{L-BTF}. One of the twisting
factors is the \textbf{jordanian twist} \cite{OGI} defined on the Borel
subalgebra $b^{2}=\{H,E\mid \lbrack H,E]=E\}$, it has the twisting element
\begin{equation*}
\mathcal{F}_{\mathcal{J}}=\exp \{H\otimes \sigma \},
\end{equation*}
where $\sigma =\ln (1+E)$. This twist generates the solution of the
Yang-Baxter equation: $\mathcal{R}=(\mathcal{F}_{\mathcal{J}})_{21}\mathcal{F%
}_{\mathcal{J}}^{-1}$. The classical $r$-matrix associated with it is $%
r=H\wedge E$.

The quantization problem for parabolic subgroups of a simple Lie group is
reduced to the problem of constructing the explicit solutions $\mathcal{F}$
(of equations (\ref{twist-equations})) whose carriers
$g_{c\left( \mathcal{F}\right) }$ are parabolic subalgebras in $g$.

In the general case the solutions of the twist equation do not form an
algebra with respect to (Hopf algebra) multiplication: 
the product of twists is not
a twist. However under certain conditions the compositions of twists
constitute the solutions of the twist equations. For infinite series of
simple Lie algebras the twists (called ``\textbf{chains}'') were found \cite
{KLO}. They were constructed as products of factors each being itself a
twisting element for the previously twisted algebra. Certain groups of
factors (''\textbf{links}'' of the chain) are twisting elements for the
initial Hopf algebra $U\left( g\right) $. The structure of these chains is
determined by the fundamental symmetry properties of the corresponding root
system and is common for all simple Lie algebras \cite{AKL}. The possibility
to compose a chain of twists is based on the existence of sequences of
regular injections
\begin{equation}
g_{m}\subset g_{m-1}\ldots \subset g_{2}\subset g_{1}=g  \label{inject}
\end{equation}
where the subalgebras $g_{k}$ belong to the same series of simple algebras
and $m$ is the \textbf{length of the chain}. The sequence (\ref{inject}) is
formed according to the following rule. Among the long roots of the root
system $\Lambda (g_{k})$ the \textbf{initial root} $\lambda _{k}^{0}$ is
chosen. Let the subspace $V_{{\lambda }_{k}^{0}}^{\bot }$ in the root space
of $g_{k}$ be orthogonal to the initial root $\lambda _{k}^{0}$ then the
root subsystem $\Lambda (g_{k+1})=\Lambda (g_{k})\bigcap V_{{\lambda }%
_{k}^{0}}^{\bot }$ defines the subalgebra $g_{k+1}$. The sets $\pi _{k}$ are
called the \textbf{constituent roots} for $\lambda _{k}^{0}$:
\begin{equation}
\begin{array}{l}
\pi _{k}=\left\{ \lambda ^{\prime },\lambda ^{\prime \prime }\mid \lambda
^{\prime }+\lambda ^{\prime \prime }=\lambda _{k}^{0};\quad \lambda ^{\prime
}+\lambda _{k}^{0},\lambda ^{\prime \prime }+\lambda _{k}^{0}\notin \Lambda
\left( g_{k}\right) \right\} \\
\pi _{k}=\pi _{k}^{\prime }\cup \pi _{k}^{\prime \prime };\qquad \pi
_{k}^{\prime }=\left\{ \lambda ^{\prime }\right\} ,\pi _{k}^{\prime \prime
}=\left\{ \lambda ^{\prime \prime }\right\} .
\end{array}
\label{c-roots}
\end{equation}
It was proved \cite{KLO,AKL} that for each sequence (\ref{inject})
corresponds the chain of twists containing $m$ links of the standard form.
For series $A_{n-1}$ consider the Cartan-Weil basis and let $E_{\lambda }\in
g$ be the vector corresponding to the root $\lambda \in \Lambda \left(
g\right) $. For each initial root $\lambda _{k}^{0}$ define the element $%
\sigma _{\lambda _{k}^{0}}=\ln \left( 1+E_{\lambda _{k}^{0}}\right) $ and
the Cartan element $H_{\lambda _{k}^{0}}$ so that $\mathrm{ad}_{H_{\lambda
_{k}^{0}}}\circ E_{\lambda _{k}^{0}}=E_{\lambda _{k}^{0}}$ and $\mathrm{ad}%
_{H_{\lambda _{k}^{0}}}\circ E_{\lambda ^{\prime }}=\frac{1}{2}E_{\lambda
^{\prime }}$. Then the product of twisting factors
\begin{equation}
\mathcal{F}_{\mathrm{link}}^{k}=\prod_{\lambda ^{\prime }\in \pi
_{k}^{\prime }}\exp \left\{ E_{\lambda ^{\prime }}\otimes E_{\lambda
_{k}^{0}-\lambda ^{\prime }}e^{-\frac{1}{2}\sigma _{\lambda
_{k}^{0}}}\right\} \cdot \exp \{H_{\lambda _{k}^{0}}\otimes \sigma _{\lambda
_{k}^{0}}\}\,  \label{link-ini}
\end{equation}
forms the link and the product of links -- the chain:
\begin{equation}
\begin{array}{l}
\mathcal{F}_{ch_{\left( 1\prec m\right) }}=\prod_{k}^{m}\mathcal{F}_{\mathrm{%
link}}^{k}
\end{array}
.  \label{chain-ini}
\end{equation}
$\mathcal{F}_{ch_{\left( 1\prec m\right) }}$\ is a solution of equations (%
\ref{twist-equations}) for the Hopf algebra $U\left( g\right) $.

Chains permit to construct twists with large carriers. For any universal
enveloping algebra $U\left( g\right) $ with simple $g$ chains of extended
jordanian twists give the possibility to construct the deformations whose
carrier subalgebras $g_{ch}\subset g$ contain maximal nilpotent subalgebras,
$g_{ch}\supset \mathbf{N}_{g}^{\pm }$ \cite{AKL}. 
Such chains are called \textbf{full}%
. Even in the case of full chains the Cartan subalgebra $\mathbf{H}\left(
g\right) $ is not always contained in $g_{ch}$. The structure of chain
is based on the solvability of its carrier \cite{AKL} and thus the only
possibility for a chain carrier $g_{ch}$ to be parabolic is to coincide with
the Borel subalgebra $\mathbf{B}^{+}(g)$.

The compliment subalgebra in $\mathbf{H}\left( g\right) $ orthogonal to all
the initial roots of the chain is denoted by $H^{\perp }$,
\begin{equation}
\mathbf{H}\left( g\right) =\left( \mathbf{H}\left( g\right) \bigcap
g_{ch}\right) \oplus H^{\perp }.  \label{Cart-subalg}
\end{equation}
It follows from the structure of the sequence (\ref{inject}) and the chain (%
\ref{chain-ini}) that the number $m$ of links in the chain is equal to the
dimension of the Cartan subalgebra in the carrier,
\begin{equation}
\dim \left( \mathbf{H}\left( g\right) \cap g_{ch}\right) =m,
\end{equation}
Obviously the number of links is not greater than the rank $r$ of $g$.
It was shown in \cite{AKL,KLS} that for $g=A_{1},B_{n}$,$C_{n}$ and $%
D_{n=2k} $ these two characteristics coincide: $m=r$.

If $\dim H^{\perp }$ is even one can always construct a nontrivial Abelian
twist $\mathcal{F}_{A}$ with the carrier $g_{A}=H^{\perp }$. It can be
easily seen that the product $\mathcal{F}_{A}\mathcal{F}_{ch}$ is again a
twist and in this case its carrier $g_{Ach}$ coincides with $\mathbf{B}%
^{+}(g)$ . This gives the desired parabolic carrier but its Levy factor is
still Abelian. The new carrier is an extension of $g_{ch}$ by an Abelian
subalgebra and thus the twists $\mathcal{F}_{A}\mathcal{F}_{ch}$ are trivial
extensions of the ordinary full chains $\mathcal{F}_{ch}$.

In this paper we shall study the special types of chains $\mathcal{F}_{%
\widehat{ch}}^{R}$ for algebras $A_{n-1}$. It will be demonstrated that for
such chains there exist such nontrivial additional factors $\mathcal{F}^{%
\widehat{J}}$ that the compositions $\mathcal{F}^{\widehat{J}}\mathcal{F}_{%
\widehat{ch}}^{R}$ are twists and their carriers $g_{\mathcal{P}}$\ are
nonminimal parabolic subalgebras of $A_{n-1}$. These parabolic carriers
cannot be reduced to an Abelian extensions of $g_{\widehat{ch}}$. They
contain $\mathbf{B}^{+}(g)$ and their intersection with $\mathbf{N}_{g}^{-}$
is nontrivial, $g_{\mathcal{P}}\cap \mathbf{N}_{g}^{-}=\emptyset $. The
corresponding parabolic twisting elements $\mathcal{F}_{\mathcal{P}}=%
\mathcal{F}^{\widehat{J}}\mathcal{F}_{\widehat{ch}}^{R}$ will be explicitly
constructed and their properties will be illustrated by examples. In the
case of $A_{2}$ the proposed construction degenerates into the so called
elementary parabolic twist presented in \cite{LSam}.

\section{What parabolic subalgebras are we interested in?}

In the previous Section it was mentioned that for algebras $B_{n}$, $C_{n}$
and  $D_{n=2k} $ there exist the twist deformations
(full chains of extended twists)
whose carrier subalgebras are minimal parabolic: $g_{ch}=\mathbf{B}^{+}(g)$.
In the case $g=A_{n-1}$ the situation is different. Here only for $sl\left(
2\right) $ the number of links coincides with the rank. In all other cases
$r$ is greater than $m$:
\begin{equation}
m=\left\{
\begin{array}{c}
\frac{n-1}{2}\qquad \text{for odd }n \\
\frac{n}{2}\qquad \text{for even }n
\end{array}
\right.  \label{links-num}
\end{equation}
Thus the full chains themselves cannot have the parabolic carriers for the
algebras $A_{n-1}$ with $n>2$.

Our aim is to construct for $U\left( A_{n-1}\right) $ the nonminimal
parabolic twists $\mathcal{F}_{\mathcal{P}}$ whose carriers $g_{\mathcal{P}}$
contain nontrivial semisimple factors. We shall search such carriers in the
form
\begin{eqnarray}
g_{\mathcal{P}} &=&s\vartriangleright \left( g_{ch}\setminus
\sum_{i=1}^{p}g^{\gamma _{i}}\right) ,  \label{semi-dir-sum1} \\
s &=&\oplus _{i=1}^{p}sl\left( 2\right) _{i},\qquad p=r_{A_{n-1}}-m,  \notag
\\
sl\left( 2\right) _{i} &=&g^{-\gamma _{i}}+\widehat{H_{i}^{\perp }}%
+g^{\gamma _{i}}.
\end{eqnarray}
Thus in the space $H^{\perp }$ there must exists a basis whose elements $%
\widehat{H_{i}^{\perp }}$ generate the Cartan subalgebras in $s=\oplus
_{s=1}^{p}sl\left( 2\right) _{s}$.

The parabolic subalgebra $g_{\mathcal{P}}$ of a simple Lie algebra $g$ with
the root system $\Lambda $ is defined by the subset $\Psi $ of the set $%
\Sigma =\left\{ \alpha _{k}\right\} $ of basic roots. The roots of $g_{%
\mathcal{P}}$ (in their standard decomposition) contain $\alpha _{j}\in \Psi
$ only with nonnegative coefficients. Let $\Phi \left( \Psi \right) \subset
\Lambda $ be the set of roots with such properties. Then the Cartan
decomposition of $g_{\mathcal{P}}$ can be written as follows,
\begin{equation}
g_{\mathcal{P}}=\mathbf{H}(g)+\sum_{\lambda \in \Phi \left( \Psi \right)
}g^{\lambda }.
\end{equation}

Let us define the sets $\Psi _{\mathcal{P}}$ so that the factor $s$ in the
semidirect sum (\ref{semi-dir-sum1}) is equivalent to the direct sum $\oplus
_{s=1}^{p}sl\left( 2\right) _{s}$. Let the basic roots of $A_{n-1}$ be
canonically expressed in the orthonormal basis $\left\{ e_{j}\right\} $
of $\mathsf{R}^{n}$%
\begin{equation}
\alpha _{k}:=e_{k}-e_{k+1}.
\end{equation}
We shall use the function
\begin{equation}
\chi \left( l\right) =\frac{1}{2}\left( n+\left( n-2l\right) \left(
-1\right) ^{l+1}\right)
\end{equation}
and fix the set $\Psi _{\mathcal{P}}$ as follows
\begin{equation}
\Psi _{\mathcal{P}}=\left\{
\begin{array}{c}
\left\{ \alpha _{\chi \left( s\right) }\mid s=1,2,\ldots ,\frac{n-1}{2}%
=p\right\} \qquad \text{for odd }n, \\
\left\{ \alpha _{\chi \left( s\right) },\alpha _{\frac{n}{2}}\mid
s=1,2,\ldots ,\frac{n-2}{2}=p\right\} \qquad \text{for even }n.
\end{array}
\right.
\end{equation}

The Cartan decomposition for the parabolic subalgebra $g_{\mathcal{P}}$ is
thus defined
\begin{equation}
g_{\mathcal{P}}=\mathbf{H}+\sum_{\lambda \in \Lambda ^{+}}g^{\lambda
}+\sum_{\gamma \in \Gamma _{\mathcal{P}}}g^{-\gamma },
\end{equation}
where
\begin{equation}
\Gamma _{\mathcal{P}}=\left\{ \alpha _{j}\mid j\left( s\right) =n-\chi
\left( s\right) ,s=1,\ldots ,p\right\}.
\end{equation}

We know how to construct full chains $\mathcal{F}_{ch}$ with the
carriers $g_{ch}$ (see Section 1). Consider the space $g_{\mathcal{P}%
}\setminus g_{ch}\equiv g_{\widehat{J}}$ . This is the space of an even
dimensional algebra that can be presented as a direct sum of Borel
subalgebras,
\begin{equation}
g_{J}=\oplus _{s}b_{s},
\end{equation}
with
\begin{equation}
b_{s}=\left\{ H_{k}^{\widehat{\perp }},g^{-\alpha _{s}}\right\} _{s=n-\chi
\left( i\right) ;\quad i=1,\ldots ,p}.
\end{equation}

In Section 5 we shall demonstrate that for the subalgebra $U_{ch}\left(
g_{\mathcal{P}}\right) $  twisted by the
full chain $\mathcal{F}_{ch}$ we can construct the product of $p$
quasi-Jordanian twists whose integral carrier contains the subalgebra $g_{J}$.
Thus the new twist
\begin{equation}
\mathcal{F}_{\mathcal{P}}=\mathcal{F}^{\widehat{J}}\mathcal{F}_{ch}
\label{parab-general}
\end{equation}
will have the parabolic carrier
\begin{equation}
g_{\mathcal{P}}=\left( \bigoplus_{s=1}^{p}sl\left( 2\right) _{i}\right)
\vartriangleright \left( g_{ch}\setminus \sum_{s=1}^{p}g^{\gamma
_{s}}\right) .
\end{equation}

\section{Full chain and dual coordinates.}

We start by constructing the special type $\mathcal{F}_{\widehat{ch}}$ of
full chain of extended jordanian twists that differs from (\ref{link-ini}),
(\ref{chain-ini}) by the choice of the Cartan elements and correspondingly by
the normalization of the extension factors. In \cite{KwL} and \cite{LO} it
was proved that such a solution of the twist equations (\ref{twist-equations})
exists:
\begin{eqnarray}
\mathcal{F}_{\widehat{ch}} &=&\prod_{l=1}^{m}\mathcal{F}_{\widehat{\mathrm{%
link}}}^{l}=  \notag \\
&=&\prod_{l=1}^{m}\left( e^{\sum_{k_{l}=l+1}^{n-l}E_{l,k_{l}}\otimes
E_{k_{l},n-l+1}e^{-\frac{\left( 1-\left( -1\right) ^{l}\right) }{2}\sigma
_{l,n-l+1}}}e^{\widehat{H_{l}}\otimes \sigma _{l,n-l+1}}\right)
\label{1-chain}
\end{eqnarray}
Remember that $m$ is the number of links in the full chain (see(\ref
{links-num})) and $\sigma _{l,n-l+1}=\ln \left( 1+E_{l,n-l+1}\right) $. The
set $\left\{ \widehat{H_{l}}\right\} $ of Cartan elements is chosen as
follows:
\begin{equation}
\begin{array}{l}
\widehat{H}_{1}=+ \frac{1}{n}\mathbf{I} -E_{n,n}; \\
\widehat{H}_{2}=- \frac{3}{n}\mathbf{I}
+E_{1,1}+E_{2,2}+E_{n,n}; \\
\ldots \\
\widehat{H}_{l}=+ \frac{2l-1}{n}\mathbf{I} - \\
\quad \quad -E_{1,1}-\ldots -E_{l-1,l-1}-E_{n-l+1,n-l+1}-\ldots
-E_{n,n};\quad l=2k-1; \\
\widehat{H}_{l}=- \frac{2l-1}{n}\mathbf{I} + \\
\quad \quad +E_{1,1}+\ldots +E_{l,l}+E_{n-l+2,n-l+2}+\ldots +E_{n,n};\quad
\quad \quad l=2k; \\
\ldots
\end{array}
\end{equation}
These elements form the Cartan subalgebra $H_{\widehat{ch}}$ and the Cartan
decomposition of the carrier $g_{\widehat{ch}}$ of the chain (\ref{1-chain})
looks like
\begin{equation*}
g_{\widehat{ch}}=H_{\widehat{ch}}\vartriangleright \left( \sum_{\lambda \in
\Lambda ^{+}}g^{\lambda }\right) .
\end{equation*}
The full chain bears the set of $m$ free parameters \cite{KLO} but here for
simplicity we shall omit them.

Applying the twist $\mathcal{F}_{\widehat{ch}}$ to the Hopf algebra $U\left(
A_{n-1}\right) $ we obtain the deformation
\begin{equation}
\mathcal{F}_{\widehat{ch}}:U\left( A_{n-1}\right) \longrightarrow U_{%
\widehat{ch}}\left( A_{n-1}\right) .
\end{equation}
According to the arguments developed in the previous Section we are to find
the solutions $\mathcal{F}^{\widehat{J}}$\ of the twist equations for $%
\Delta _{\widehat{ch}}$. The carrier of such solutions is to contain the
direct sum of 2-dimensional Borel subalgebras: $g_{J}=\oplus _{s}b_{s}$.

Usually the Hopf structure of the deformed algebras such as $U_{\widehat{ch}%
}\left( A_{n-1}\right) $ is explicitly described in terms of matrix
coordinate functions $E_{ij}$ and the normalized Cartan generators
\begin{equation}
H_{ij}:=\frac{1}{2}\left( E_{ii}-E_{jj}\right) .
\end{equation}
Despite the property of selfduality of twisted algebras $U_{\mathcal{F}%
}\left( A_{n-1}\right) $ the coordinates $E_{ij}$ are not convenient for the
description of the costructure. For example, the Jordanian deformation $U_{%
\mathcal{J}}\left( A_{1}\right) $ is always written in terms of $\sigma
_{12}$ but not $E_{12}$. For the deformed Hopf algebra $U_{\mathcal{F}%
}\left( A_{n-1}\right) $ the most appropriate is the basis consisting of the
dual group coordinates. The dual group basis for twisted universal
enveloping algebras was studied in \cite{L}, there the explicit procedure
was proposed to find the relations between the algebraic coordinates (of the
type $E_{ij}$) and the dual coordinates $\widehat{E}_{ij}$. As we shall see
below the Borel subalgebras $b_{s}$ in $g_{P}\subset U_{\widehat{ch}}\left(
A_{n-1}\right) $ with the necessary Hopf structure could be defined only in
the deformed (dual group) basis.

As far as the general form of the link factor $\mathcal{F}_{\widehat{%
\mathrm{link}}}$ of the proposed twists (\ref{parab-general}) is known (up
to the Abelian ''rotation'' that will be fixed later) we shall focus our
attention on the dual coordinates for the subalgebra $U_{\widehat{ch}}\left(
g_{J}\right) $.

Let us remind \cite{KLO} that the full chain does not change the costructure
on the subspace $H^{\perp }$ (see (\ref{Cart-subalg})). This means that on
this space the dual group coordinates coincide with the algebraic ones. We
denote them by $H_{s}^{\widehat{\perp }}$ and choose the following set as a
basis in $H^{\perp }$:
\begin{eqnarray}
H_{s}^{\widehat{\perp }} &=&\left( -1\right) ^{s}\left( -\frac{2s}{n}%
\mathbf{I}+\sum_{u=1}^{s}\left( E_{u,u}+E_{n-u+1,n-u+1}\right)
\right) ; \\
s &=&1,\ldots ,p.  \notag
\end{eqnarray}

For the generators of the subspaces $g^{-\alpha _{j}},$ $j\left( s\right)
=n-\chi \left( s\right) ,s=1,\ldots ,p,\quad $ we apply the algorithms
developed in \cite{L} and obtain the following expressions for the dual
coordinates $\widehat{E}_{s}$. Notice that the forms of coordinates are
different for odd and even indices, in what follows we shall call them 
"odd"
and "even external coordinates" (to distinguish them from the "internal"
coordinates defined on $B^{+}$);
\begin{eqnarray}
\widehat{E}_{s}\! &=&\!\left\{
\begin{array}{c}
E_{j\left( s\right) +1,j\left( s\right) }+E_{\chi \left( s\right) ,\chi
\left( s\right) +1}e^{\sigma _{j\left( s\right) +1,\chi \left( s\right)
}-\sigma _{j\left( s\right) ,\chi \left( s\right) +1}};\quad s=2k-1; \\
E_{j\left( s\right) +1,j\left( s\right) }\!+\!\left( \widehat{H}_{\chi
\left( s\right) +1}\!-\!\widehat{H}_{\chi \left( s\right) }\right) \!E_{\chi
\left( s\right) ,j\left( s\right) }\!+\!E_{\chi \left( s\right) ,\chi \left(
s\right) +1};~s=2k.
\end{array}
\right.  \label{coordinates-prelim} \\
s &=&1,\ldots ,\left\{
\begin{array}{c}
p-1\quad \text{for odd }n \\
p\quad \text{for even }n
\end{array}
\right. .  \notag
\end{eqnarray}
When $n$ is odd the form of the ''last'' external coordinate $\widehat{E}%
_{p} $ degenerates:
\begin{equation}
\widehat{E}_{p}=\left\{
\begin{array}{c}
E_{j\left( p\right) +1,j\left( p\right) };\quad \text{for odd }n\text{ and
odd }p \\
E_{j\left( p\right) +1,j\left( p\right) }-\widehat{H}_{\chi \left( p\right)
}E_{\chi \left( p\right) ,j\left( p\right) };\quad \text{for odd }n\text{
and even }p
\end{array}
\right. .  \notag
\end{equation}

As far as
\begin{eqnarray}
\quad \text{for odd }s\quad j\left( s\right) &=&s,  \notag \\
\quad \text{for even }s\quad \chi \left( s\right) &=&s,
\end{eqnarray}
these expressions can be rewritten in a more compact way:
\begin{eqnarray}
\widehat{E}_{s} &=&\left\{
\begin{array}{c}
E_{s+1,s}+E_{n-s,n-s+1}e^{\sigma _{s+1,n-s}-\sigma _{s,n-s+1}};\quad s=2k-1;
\\
E_{n-s+1,n-s}+\left( \widehat{H}_{s+1}-\widehat{H}_{s}\right)
E_{s,n-s}+E_{s,s+1};\quad \quad s=2k;
\end{array}
\right.  \label{coordinates-short} \\
s &=&1,\ldots ,\left\{
\begin{array}{c}
p-1\quad \text{for odd }n \\
p\quad \text{for even }n
\end{array}
\right. .  \notag
\end{eqnarray}

\begin{equation}
\widehat{E}_{p}=\left\{
\begin{array}{c}
E_{p+1,p};\quad \text{for odd }n\text{ and odd }p \\
E_{n-p+1,n-p}-\widehat{H}_{p}E_{p,n-p};\quad \text{for odd }n\text{ and even
}p
\end{array}
\right. .  \notag
\end{equation}

\section{Properties of external coordinates}

\begin{conjecture}
Let $\mathcal{L}\left( \widehat{E}_{s}\right) $ be the space spanned by the
set of external coordinates. Then $H^{\perp }\oplus \mathcal{L}\left(
\widehat{E}_{s}\right) $ is the space of the direct sum $g_{J}=\oplus
_{s}b_{s}$ . Each 2-dimensional Borel subalgebra $b_{s}$ is generated by the
corresponding pair $\left\{ H_{s}^{\widehat{\perp }},\widehat{E}_{s}\right\}
$.
\end{conjecture}

\begin{proof}
This can be checked by direct computations. All the external coordinates
commute,
\begin{equation}
\left[ \widehat{E}_{s},\widehat{E}_{t}\right] =0;\quad s,t=1,\ldots ,p.
\end{equation}
Notice that the subsets of odd external coordinates and of even coordinates
are trivially commutative. This obviously follows from the form of their
expressions in terms of the elements $E_{i,j}$ (see (\ref{coordinates-short}%
)): the ranges of indices $i,j$ in different $\widehat{E}_{s}$ in these two
subsets do not intersect. Thus only the vanishing of odd-even commutators is
to be checked. It also can be verified that the Cartan elements $H_{s}^{%
\widehat{\perp }}$ act on the external coordinates as follows
\begin{equation}
\left[ H_{s}^{\widehat{\perp }},\widehat{E}_{t}\right] =\delta _{s,t}%
\widehat{E}_{t}.  \label{H-E-comm-2}
\end{equation}
\end{proof}

There exists an important relation connecting the Cartan generators 
$\widehat{%
H_{s}}$ in the Jordanian factors $e^{\widehat{H_{l}}\otimes \sigma
_{l,n-l+1}}$ of the chain $\mathcal{F}_{\widehat{ch}}$ and the basic
elements $H_{s}^{\widehat{\perp }}$:
\begin{equation}
-2H_{s}^{\widehat{\perp }}=\left\{
\begin{array}{c}
2H_{s,s+1}+\widehat{H}_{s+1}-\widehat{H}_{s},\quad s=2k-1, \\
2H_{n-s,n-s+1}+\widehat{H}_{s+1}-\widehat{H}_{s}\quad s=2k.
\end{array}
\right.  \label{cartans-interrel}
\end{equation}

Now consider the costructure induced by the twist deformation $\mathcal{F}_{%
\widehat{ch}}$ (see (\ref{1-chain})) on the subalgebra generated by $$
\left\{ H_{s}^{\widehat{\perp }},\widehat{E}_{t},\sigma
_{l,n-l+1}\mid s,t=1,\ldots ,p;l=1,\ldots ,m\right\}. $$

\begin{conjecture}
The only nontrivial cocommutator in $U_{\widehat{ch}}\left( A_{n-1}\right) $
with values in $\mathcal{L}\left( \widehat{E}_{s}\right) $ is a map $\left(
g^{\alpha _{j\left( s\right) }}\right) ^{\#}\otimes \left( H_{s}^{\widehat{%
\perp }}\right) ^{\#}\longrightarrow \left( \widehat{E}_{s}\right) ^{\#}$.
The corresponding coproducts are
\begin{multline}
\Delta _{\widehat{ch}}\left( \widehat{E}_{s}\right) = \\
=\left\{
\begin{array}{c}
\widehat{E}_{s}\otimes 1+1\otimes \widehat{E}_{s}-2H_{s}^{\widehat{\perp }%
}\otimes E_{s+1,n-s+1}e^{-\sigma _{s,n-s+1}};s=2k-1; \\
\widehat{E}_{s}\otimes 1+1\otimes \widehat{E}_{s}-E_{s,n-s}\otimes 2H_{s}^{%
\widehat{\perp }};\quad s=2k.
\end{array}
\right.   \label{ext-coprs}
\end{multline}
\end{conjecture}

\begin{proof}
From (\ref{1-chain}) it can be observed that for any odd indicae $s$ the
element $E_{s+1,n-s+1}e^{-\sigma _{s,n-s+1}}$ is dual to $E_{s,s+1}$. The
functional corresponding to $E_{s+1,n-s+1}e^{-\sigma _{s,n-s+1}}$ acts
inversely with respect to $E_{s,s+1}$ and thus can shift the elements from $%
H^{\perp }$ to $\mathcal{L}\left( \widehat{E}_{s}\right) $. Correspondingly
for even $s$ the dual for $E_{s,n-s}$ is $E_{n-s,n-s+1}$ with the similar
shifting properties. Notice that the appearance of the factor $e^{-\sigma
_{s,n-s+1}}$ in the odd case is in complete correspondence with the form of
odd and even links in the chain $\mathcal{F}_{\widehat{ch}}$.

Let us check the relation (\ref{ext-coprs}).

Analyzing the explicit form of the external coordinates $\widehat{E}_{s}$
one can see that each of them is defined on the subalgebra

\begin{equation}
M_{s}=\left( sl\left( 2\right) _{s,s+1}\oplus sl\left( 2\right)
_{n-s,n-s+1}\right) \vdash J_{s}
\end{equation}
with $J_{s}$ being an ideal corresponding to the roots: $%
e_{s}-e_{n-s}$, $e_{s}-e_{n-s+1}$, $e_{s+1}-e_{n-s}$\ and $e_{s}-e_{n-s+1}$.
Under the action of first $s-1$ links of $\mathcal{F}_{\widehat{ch}}$ the
coproducts in $M_{s}$ remain unchanged (primitive on the generators) due to
the ''matreshka'' effect \cite{KLO}. On the other hand the links $\mathcal{F}%
_{\widehat{\mathrm{link}}}^{l}$ with $l>s+1$ are defined on the subalgebras
\begin{equation}
sl^{\vee }(n-4s)
\end{equation}
with the roots $e_{i}-e_{j}$; $i,j=s+2,\ldots ,n-s-1$. In $sl(n)$ these
subalgebras, $M_{s}$ and $sl^{\vee }(n-4s)$, form a direct sum:
\begin{equation}
sl(n)\supset M_{s}\oplus sl^{\vee }(n-4s).
\end{equation}
It follows that the links $\mathcal{F}_{l>s+1}$ also cannot deform the
costructure of the coordinate $\widehat{E}_{s}$. Thus only two links are to
be considered for each $\widehat{E}_{s}$, the $s$-th and
the $\left( s+1\right)$-th:

\begin{equation}
\mathcal{F}_{s}=\exp \left( \sum_{k_{s}=s+1}^{n-s}E_{s,k_{s}}\otimes
E_{k_{s},n-s+1}e^{-\sigma _{s,n-s+1}}\right) \exp \left( \widehat{H_{s}}%
\otimes \sigma _{s,n-s+1}\right)
\end{equation}
\begin{equation}
\mathcal{F}_{s+1}=\exp \left(
\sum_{k_{s+1}=s+2}^{n-s-1}E_{s+1,k_{s+1}}\otimes E_{k_{s+1},n-s}\right) \exp
\left( \widehat{H_{s+1}}\otimes \sigma _{s+1,n-s}\right)
\end{equation}

Let us apply the corresponding adjoint operators to the initial primitive
coproducts $\Delta ^{\left( 0\right) }\left( E_{s+1,s}\right) $ and obtain
for them the twisted expressions:
\begin{equation}
\begin{array}{l}
\Delta _{s,s+1}\left( E_{s+1,s}\right) =\mathcal{F}_{s+1}\circ \mathcal{F}%
_{s}\circ \left( E_{s+1,s}\otimes 1+1\otimes E_{s+1,s}\right) \rule{0mm}{6mm}%
= \\
\rule{0mm}{6mm}\quad =e^{\mathrm{ad}%
\sum_{k_{s+1}=s+2}^{n-s-1}E_{s+1,k_{s+1}}\otimes E_{k_{s+1},n-s}}e^{\mathrm{%
ad}\widehat{H_{s+1}}\otimes \sigma _{s+1,n-s}}\circ  \\
\rule{0mm}{6mm}\quad \quad \circ e^{\mathrm{ad}%
\sum_{k_{s}=s+1}^{n-s}E_{s,k_{s}}\otimes E_{k_{s},n-s+1}e^{-\sigma
_{s,n-s+1}}}e^{\mathrm{ad}\left( \widehat{H_{s}}\otimes \sigma
_{s,n-s+1}\right) } \\
\rule{0mm}{6mm}\quad \quad \circ \left( E_{s+1,s}\otimes 1+1\otimes
E_{s+1,s}\right) = \\
\\
\quad =E_{s+1,s}\otimes 1+1\otimes E_{s+1,s}+ \\
\rule{0mm}{6mm}\quad \quad +\left( 2H_{s,s+1}+\widehat{H_{s+1}}-\widehat{%
H_{s}}\right) \otimes E_{s+1,n-s+1}e^{-\sigma _{s,n-s+1}}+ \\
\rule{0mm}{6mm}\quad \quad
-\sum_{k_{s+1}=s+2}^{n-s-1}E_{s+1,k_{s+1}}e^{\sigma _{s+1,n-s}}\otimes
E_{k_{s+1},n-s+1}e^{\sigma _{s+1,n-s}-\sigma _{s,n-s+1}} \\
\rule{0mm}{6mm}\quad \quad
+\sum_{k_{s+1}=s+2}^{n-s-1}E_{s+1,k_{s+1}}e^{\sigma _{s+1,n-s}}\otimes
E_{k_{s+1},n-s}E_{s+1,n-s+1}e^{-\sigma _{s,n-s+1}} \\
\rule{0mm}{6mm}\quad \quad -e^{\sigma _{s+1,n-s}}\otimes
E_{n-s,n-s+1}e^{\sigma _{s+1,n-s}-\sigma _{s,n-s+1}} \\
\rule{0mm}{6mm}\quad \quad +1\otimes E_{n-s,n-s+1}e^{\sigma
_{s+1,n-s}-\sigma _{s,n-s+1}} \\
\rule{0mm}{6mm}\quad \quad -\widehat{H_{s+1}}e^{\sigma _{s+1,n-s}}\otimes
E_{s+1,n-s+1}e^{-\sigma _{s,n-s+1}};
\end{array}
\end{equation}
The result of the twist deformation by $\mathcal{F}_{s+1}\circ \mathcal{F}%
_{s}$ for the exponential factor $e^{\sigma _{s+1,n-s}-\sigma _{s,n-s+1}}$
in the second summand of $\widehat{E}_{s}$ is well known \cite{KLM} -- this
factor is group-like:
\begin{equation}
\Delta _{s,s+1}\left( e^{\sigma _{s+1,n-s}-\sigma _{s,n-s+1}}\right)
=e^{\sigma _{s+1,n-s}-\sigma _{s,n-s+1}}\otimes e^{\sigma _{s+1,n-s}-\sigma
_{s,n-s+1}};
\end{equation}
Thus only the coproduct $\Delta _{s,s+1}\left( E_{n-s,n-s+1}\right) $
remains to be constructed:
\begin{equation}
\begin{array}{l}
\Delta _{s,s+1}\left( E_{n-s,n-s+1}\right) =\mathcal{F}_{s+1}\circ \mathcal{F%
}_{s}\circ \left( E_{n-s,n-s+1}\otimes 1+1\otimes E_{n-s,n-s+1}\right) =\rule%
{0mm}{6mm} \\
\rule{0mm}{6mm}\quad =E_{n-s,n-s+1}\otimes e^{\sigma _{s,n-s+1}-\sigma
_{s+1,n-s}}+e^{\sigma _{s,n-s+1}}\otimes E_{n-s,n-s+1}+ \\
\rule{0mm}{6mm}\quad \quad
+\sum_{k_{s+1}=s+2}^{n-s-1}E_{s+1,k_{s+1}}e^{\sigma _{s,n-s+1}}\otimes
E_{k_{s+1},n-s+1}+ \\
\rule{0mm}{6mm}\quad \quad +\widehat{H_{s+1}}e^{\sigma _{s,n-s+1}}\otimes
E_{s+1,n-s+1}e^{-\sigma _{s+1,n-s}}- \\
\rule{0mm}{6mm}\quad \quad
-\sum_{k_{s+1}=s+2}^{n-s-1}E_{s+1,k_{s+1}}e^{\sigma _{s,n-s+1}}\otimes
E_{k_{s+1},n-s}E_{s+1,n-s+1}e^{-\sigma _{s+1,n-s}};
\end{array}
\end{equation}
As we have expected the cumbersome parts of these coproducts cancel in the
integral expression $\Delta _{s,s+1}\left( \widehat{E}_{s}\right) $ and we
obtain for $s=2k-1$
\begin{eqnarray}
\Delta _{\widehat{ch}}\left( \widehat{E}_{s}\right)  &=&\Delta
_{s,s+1}\left( \widehat{E}_{s}\right) =  \notag \\
-&=&\Delta _{s,s+1}\left( E_{s+1,s}+E_{n-s,n-s+1}e^{\sigma _{s+1,n-s}-\sigma
_{s,n-s+1}}\right) =  \notag \\
&=&\widehat{E}_{s}\otimes 1+1\otimes \widehat{E}_{s}+  \notag \\
&&+\left( 2H_{s,s+1}+\widehat{H_{s+1}}-\widehat{H_{s}}\right) \otimes
E_{s+1,n-s+1}e^{-\sigma _{s,n-s+1}}.
\end{eqnarray}

For even coordinates $\widehat{E}_{s=2k}=E_{n-s+1,n-s}+\left( \widehat{H}%
_{s+1}-\widehat{H}_{s}\right) E_{s,n-s}+E_{s,s+1}$ the nontrivial part of
the deformation is performed by the links
\begin{eqnarray}
\mathcal{F}_{s+1}\mathcal{F}_{s}
&=&e^{\sum_{k_{s+1}=s+2}^{n-s-1}E_{s+1,k_{s+1}}\otimes
E_{k_{s+1},n-s}e^{-\sigma _{s+1,n-s}}}e^{\widehat{H_{s+1}}\otimes \sigma
_{s+1,n-s}}\cdot   \notag \\
&&\cdot e^{\sum_{k_{s}=s+1}^{n-s}E_{s,k_{s}}\otimes
E_{k_{s},n-s+1}}e^{\left( \widehat{H_{s}}\otimes \sigma _{s,n-s+1}\right) },
\\
s &=&2k.  \notag
\end{eqnarray}
Here is the result of the chain
deformation for even $s$ $\left( \widehat{E}_{s=2k}\right) $
\begin{multline}
\Delta _{\widehat{ch}}\left( \widehat{E}_{s}\right) =\Delta
_{s,s+1}\left( \widehat{E}_{s}\right) =  \notag \\
=\Delta _{s,s+1}\left( E_{n-s+1,n-s}+\left( \widehat{H}_{s+1}-\widehat{H}%
_{s}\right) E_{s,n-s}+E_{s,s+1}\right) =  \notag \\
=\widehat{E}_{s}\otimes 1+1\otimes \widehat{E}_{s}+E_{s,n-s}\otimes \left(
2H_{n-s,n-s+1}+\widehat{H}_{s+1}-\widehat{H}_{s}\right) ;
\end{multline}
Taking into account the relations (\ref{cartans-interrel}) we can rewrite
the obtained coproducts:
\begin{equation}
\Delta _{\widehat{ch}}\left( \widehat{E}_{s}\right) =\left\{
\begin{array}{c}
\widehat{E}_{s}\otimes 1+1\otimes \widehat{E}_{s}-2H_{s}^{\widehat{\perp }%
}\otimes E_{s+1,n-s+1}e^{-\sigma _{s,n-s+1}},\quad s=2k-1; \\
\widehat{E}_{s}\otimes 1+1\otimes \widehat{E}_{s}-E_{s,n-s}\otimes 2H_{s}^{%
\widehat{\perp }},\qquad \qquad s=2k.
\end{array}
\right.   \label{Deltas-es-chain}
\end{equation}
\end{proof}

\section{Rotations for the chains of twists.}

Now let us demonstrate that the twisted coproducts $\Delta _{\widehat{ch}%
}\left( \widehat{E}_{s}\right) $ can be reduced further with the help of
additional Abelian twists \cite{RESH}.

We have found out that the only nontrivial terms in $\Delta _{\widehat{ch}%
}\left( \widehat{E}_{s}\right) $ are due to the coadjoint action of $\left(
E_{s+1,n-s+1}e^{-\sigma _{s,n-s+1}}\right) ^{\#}$ and $\left(
E_{n-s,n-s+1}\right) ^{\#}$. The corresponding coadjoint maps (dual to the
adjoint action of $\mathrm{ad}\left( E_{s,s+1}\right) $ and $\mathrm{ad}%
\left( E_{s,n-s}\right) $) depend on the dualization $\left( \widehat{H_{s}}%
\right) ^{\#}\sim \sigma _{s,n-s+1}$. In the chain $\mathcal{F}_{\widehat{ch}%
}$ this dualization depends on the form of the Jordanian factors $e^{%
\widehat{H_{l}}\otimes \sigma _{l,n-l+1}}$. When the carrier of the chain
does not cover the Cartan subalgebra the relation $\left( \widehat{H_{s}}%
\right) ^{\#}\sim \sigma _{s,n-s+1}$ can be transformed by applying the
Abelian twists of the type $e^{\gamma _{s}H_{s}^{\widehat{\perp }}\otimes
\sigma _{s,n-s+1}}$. Such a transformation takes place for example when a
canonical chain is changed into a peripheric one \cite{KwL}. Varying the
parameters $\gamma _{s}$ one can alter the operators $\mathrm{coad}\left(
E_{s+1,n-s+1}e^{-\sigma _{s,n-s+1}}\right) ^{\#}$ and $\mathrm{coad}\left(
E_{n-s,n-s+1}\right) ^{\#}$ and under some additional conditions trivialize
their action on the
space $\left( H^{\perp }\right) ^{\#}$. Thus studying the dual algebra
action we see that there can exist a possibility to annulate the nontrivial
terms in $\Delta _{\widehat{ch}}\left( \widehat{E}_{s}\right) $ with the
help of Abelian twists.

When the Abelian twists are of the pure form $e^{\gamma _{s}H_{s}^{\widehat{%
\perp }}\otimes \sigma _{s,n-s+1}}$ they can be incorporated into the
initial Jordanian factors of the chain and their application reduces to the
transition from one type of the chain to the other \cite{KwL}. In the
present case this will be impossible because the necessary Abelian twists
have more general configuration.

Let us introduce the set
\begin{equation}
C_{s}=e^{2\sum_{i=1}^{s}\sigma _{i,n-i+1}}
\end{equation}
with the commutation properties

\begin{equation}
\left( -1\right) ^{s}\left[ \ln \left( C_{s}\right) ,\widehat{E}_{s}\right]
=2\delta _{s,t}\cdot \left\{
\begin{array}{c}
E_{s+1,n-s+1}e^{-\sigma _{s,n-s+1}};\quad \text{for odd }s \\
E_{s,n-s};\quad \text{for even }s
\end{array}
\right.   \label{E-K-commutation}
\end{equation}
\begin{eqnarray}
\left[ \ln \left( C_{t}\right) ,E_{s+1,n-s+1}e^{-\sigma _{s,n-s+1}}\right]
&=&0;  \label{C-E-comm-1} \\
\left[ \ln \left( C_{t}\right) ,E_{s,n-s}\right]  &=&0;  \notag \\
\left[ C_{s},C_{t}\right]  &=&0;  \label{K-commutations}
\end{eqnarray}
and consider the rotation twists:
\begin{eqnarray}
\mathcal{F}_{s}^{R} &=&e^{-H_{s}^{\widehat{\perp }}\otimes \ln \left(
C_{s}\right) };\quad \text{for odd }s, \\
\mathcal{F}_{s}^{R} &=&e^{\ln \left( C_{s}\right) \otimes H_{s}^{\widehat{%
\perp }}};\quad \text{for even }s.
\end{eqnarray}
Due to the relations (\ref{K-commutations}) and the obvious property
\begin{equation}
\left[ H_{s}^{\widehat{\perp }},\ln \left( C_{s}\right) \right] =0
\label{H-C-comm}
\end{equation}
the factors $\mathcal{F}_{s}^{R}$ commute with each other. According to the
relation (\ref{E-K-commutation}) we also have
\begin{equation}
\left[ \mathcal{F}_{s}^{R},\Delta _{\widehat{ch}}\left( \widehat{E}%
_{t}\right) \right] _{s\neq t}=0.
\end{equation}
It follows that the factors $\mathcal{F}_{s}^{R}$ can be applied to the Hopf
algebra $U_{\widehat{ch}}\left( A_{n-1}\right) $ and the rotation twist can
be composed as a product:
\begin{equation}
\mathcal{F}^{R}=\prod_{s=1}^{p}\mathcal{F}_{s}^{R}=\exp\left(\sum_{\mathrm{even}%
\; s}\ln \left( C_{s}\right) \otimes H_{s}^{\widehat{\perp }}-\sum_{%
\mathrm{odd}\; s}H_{s}^{\widehat{\perp }}\otimes \ln \left( C_{s}\right)
\right).
\end{equation}
Each factor $\mathcal{F}_{s}^{R}$ in $\mathcal{F}^{R}$ acts nontrivially
only on the coproduct of the corresponding external coordinate $\Delta _{%
\widehat{ch}}\left( \widehat{E}_{s}\right) $.

Taking into account the relations (\ref{H-E-comm-2}) (\ref{C-E-comm-1}) and (%
\ref{H-C-comm}) we can construct the deformation of the coproducts $\Delta _{%
\widehat{ch}}\left( \widehat{E}_{s}\right) $ induced by the rotation $%
\mathcal{F}^{R}$:
\begin{equation}
\Delta _{\widehat{ch}}^{R}\left( \widehat{E}_{s}\right) =\mathcal{F}%
^{R}\circ \Delta _{\widehat{ch}}\left( \widehat{E}_{s}\right) =\left\{
\begin{array}{c}
\widehat{E}_{s}\otimes C_{s}^{-1}+1\otimes \widehat{E}_{s},\text{\quad }%
s=2k-1, \\
\widehat{E}_{s}\otimes 1+C_{s}\otimes \widehat{E}_{s},\quad s=2k.
\end{array}
\right. ;  \label{q-primitive-coord}
\end{equation}
For the coproducts $\Delta _{\widehat{ch}}\left( H_{s}^{\widehat{\perp }%
}\right) $ the result of twisting by $\mathcal{F}^{R}$ is trivial:
\begin{equation}
\Delta _{\widehat{ch}}^{R}\left( H_{s}^{\widehat{\perp }}\right) =\mathcal{F}%
^{R}\circ \Delta ^{\left( 0\right) }\left( H_{s}^{\widehat{\perp }}\right)
=\Delta ^{\left( 0\right) }\left( H_{s}^{\widehat{\perp }}\right) .
\label{primitive-H}
\end{equation}
Thus the rotated chain
\begin{equation}
\mathcal{F}_{\widehat{ch}}^{R}=\mathcal{F}^{R}\mathcal{F}_{\widehat{ch}}
\end{equation}
performs the Hopf algebra deformation
\begin{equation}
\mathcal{F}_{\widehat{ch}}^{R}:U\left( A_{n-1}\right) \longrightarrow U_{%
\widehat{ch}}^{R}\left( A_{n-1}\right).
\end{equation}

\section{The quasi-Jordanian factors and the parabolic twists}

The twisted algebra $U_{\widehat{ch}}^{R}\left( A_{n-1}\right) $ contains
the subalgebra generated by the primitive $H_{s}^{\widehat{\perp }}$,
group-like $C_{s}$ and quasiprimitive external coordinates $\widehat{E}_{s}$.
At this point it is convenient to redefine the external coordinates. Let us
introduce
\begin{equation}
D_{s}=\left\{
\begin{array}{c}
\widehat{E}_{s}C_{s},\quad s=2k-1, \\
\widehat{E}_{s},\quad s=2k,
\end{array}
\right.
\end{equation}
and consider in $U_{\widehat{ch}}^{R}\left( A_{n-1}\right) $ the subalgebras
$\frak{B}_{s}$ generated by the triples $\left\{ H_{s}^{\widehat{\perp }%
},\right. $ $\left. D_{s},C_{s}\right\} $ and presented by the relations
\begin{eqnarray}
\left[ H_{s}^{\widehat{\perp }},D_{t}\right]  &=&\delta _{st}D_{t};\qquad
\left[ H_{s}^{\widehat{\perp }},C_{t}\right] =0;\qquad s,t=1,\ldots ,p;
\label{b2-rel-1} \\
\Delta _{\widehat{ch}}^{R}\left( H_{s}^{\widehat{\perp }}\right)  &=&\Delta
^{\left( 0\right) }\left( H_{s}^{\widehat{\perp }}\right) =H_{s}^{\widehat{%
\perp }}\otimes 1+1\otimes H_{s}^{\widehat{\perp }};  \label{b2-rel-2} \\
\Delta _{\widehat{ch}}^{R}\left( D_{s}\right)  &=&\mathcal{F}^{R}\circ
\Delta _{\widehat{ch}}\left( D_{s}\right) =D_{s}\otimes 1+C_{s}\otimes D_{s};
\label{b2-rel-3} \\
\Delta _{\widehat{ch}}^{R}\left( C_{s}\right)  &=&C_{s}\otimes C_{s}.
\label{b2-rel-4}
\end{eqnarray}
Notice that for $s\neq t$ the generators $C_{s}$ and $D_{t}$ commute (as a
consequence of (\ref{E-K-commutation})). It follows that
\begin{equation}
\begin{array}{l}
\left[ D_{s},D_{t}\right] =0; \\
\left[ D_{s}+C_{s},D_{t}+C_{t}\right] =0;
\end{array}
\quad \mathrm{for\ any}\quad s,t.  \label{D,C+D-commutativity}
\end{equation}

The factor $C_{s}$ in (\ref{b2-rel-3}) does not allow the canonical
Jordanian twist to be a solution to the twist equations for $\Delta _{%
\widehat{ch}}^{R}$ . To overcome this difficulty the following lemma will be
proved.

\begin{lemma}
Let $\frak{B}$ be a Hopf algebra generated
by the elements $H,C,D $ with the properties:

\begin{itemize}
\item
\begin{equation}
\Delta \left( H\right) =\Delta ^{\left( 0\right) }\left( H\right)
\label{condition1}
\end{equation}

\item
\begin{equation}
\Delta \left( C+D\right) =D\otimes 1+C\otimes \left( C+D\right) .
\label{condition2}
\end{equation}

\item
\begin{equation}
\left[ H,D\right] =D,\quad \left[ H,C\right] =0,  \label{condition3}
\end{equation}

Then
\begin{equation*}
\mathcal{F}=e^{H\otimes \ln \left( C+D\right) }
\end{equation*}
is a solution to the twist equations for $\frak{B}$ .
\end{itemize}
\end{lemma}

\begin{proof}
Consider
\begin{eqnarray*}
\Delta _{\mathcal{F}}\left( C+D\right)  &=&\mathcal{F}\circ \left( \Delta
\left( C+D\right) \right) = \\
\mathcal{F}\circ \left( D\otimes 1+C\otimes \left( C+D\right) \right)
&=&D\otimes \left( C+D\right) +C\otimes \left( C+D\right) = \\
&=&\left( C+D\right) \otimes \left( C+D\right) ;
\end{eqnarray*}
This means that for
\begin{equation*}
\omega =\ln \left( C+D\right)
\end{equation*}
we have
\begin{equation}
\Delta _{\mathcal{F}}\left( \omega \right) =\Delta ^{\left( 0\right) }\left(
\omega \right)
\end{equation}
and taking into account the primitivity of $\Delta ^{\left( 0\right) }\left(
H\right) $ the factorizable twist equations
\begin{eqnarray*}
\left( \mathrm{id}\otimes \Delta _{\mathcal{F}}\right) \mathcal{F} &=&\left(
\mathcal{F}\right) _{12}\left( \mathcal{F}\right) _{13}; \\
\left( \Delta \otimes \mathrm{id}\right) \mathcal{F} &=&\left( \mathcal{F}%
\right) _{13}\left( \mathcal{F}\right) _{23}.
\end{eqnarray*}
are valid for $\mathcal{F}$.
\end{proof}

\begin{remark}
In what follows we shall consider the group-like $C$. But in the general
case the class of algebras that meet the conditions of the Lemma has no
restrictions on $\Delta \left( C\right) $ and the composition $\left[ C,D%
\right] $.
\end{remark}

\begin{remark}
For a particular case $C=1$ the subalgebra $\frak{B}$ is reduced to $b^{2}$
and the element $\mathcal{F}$ becomes the ordinary Jordanian twist.
\end{remark}

For any $s=1,\ldots ,p$ the triple $\left\{ H_{s}^{\widehat{\perp }%
},D_{s},C_{s}\right\} $ comply with the conditions (\ref{condition1}-\ref
{condition3}). Thus for any $s$ we have a subalgebra $\frak{B}_{s}$
generated by $\left\{ H_{s}^{\widehat{\perp }},D_{s},C_{s}\right\} $ and can
construct the twisting element
\begin{equation}
\mathcal{F}_{s}^{\widehat{J}}=e^{H_{s}^{\widehat{\perp }}\otimes \ln \left(
C_{s}+D_{s}\right) }=e^{H_{s}^{\widehat{\perp }}\otimes \omega _{s}}
\end{equation}
that can be applied to $U_{\widehat{ch}}^{R}\left( A_{n-1}\right) $
\begin{equation}
\mathcal{F}_{s}^{\widehat{J}}:U_{\widehat{ch}}^{R}\left( A_{n-1}\right)
\longrightarrow U_{\widehat{ch}}^{\widehat{J}s}\left( A_{n-1}\right) .
\end{equation}
In the twisted subalgebra $\mathcal{F}_{s}^{\widehat{J}}\circ \frak{B}_{s}=%
\frak{B}_{s}^{\widehat{J}}$ the element $\omega _{s}$ becomes primitive
\begin{equation}
\Delta _{\widehat{ch}}^{\widehat{J}s}\left( \omega _{s}\right) =\omega
_{s}\otimes 1+1\otimes \omega _{s}.
\end{equation}
The Hopf algebra $U_{\widehat{ch}}^{\widehat{J}s}\left( A_{n-1}\right) $ can
be considered as the deformation of $U\left( A_{n-1}\right) $ by the
factorizable twist
\begin{equation}
\mathcal{F}_{s}^{\widehat{J}}\mathcal{F}_{\widehat{ch}}^{R}=\mathcal{F}_{s}^{%
\widehat{J}}\mathcal{F}^{R}\mathcal{F}_{\widehat{ch}}
\end{equation}
with the carrier
\begin{equation}
g_{ch}^{\widehat{J}s}=\left( sl\left( 2\right) \right) _{s}\vartriangleright
\left( g_{ch}\setminus g^{\gamma _{s}}\right) .
\end{equation}

The question arises whether the factors $\mathcal{F}_{t}^{\widehat{J}}$ can
be applied to $U_{\widehat{ch}}^{\widehat{J}s}\left( A_{n-1}\right) $ with $%
s\neq t$ and finally whether the same can be done with the product of all
the factors $\mathcal{F}_{s}^{\widehat{J}}$ to obtain finally the twist with
the parabolic carrier $g_{\mathcal{P}}$.

Consider the action of $\mathcal{F}_{s}^{\widehat{J}}$ on the subalgebra $%
\frak{B}_{t}$ with $s\neq t$. The invariance of the coproducts (\ref
{condition1}) and (\ref{condition2}) with respect to the twist $\mathcal{F}%
_{s}^{\widehat{J}}$ is needed. The coproducts of $H_{t}^{\perp }\mid _{t\neq
s}$ are invariant due to the relations (\ref{b2-rel-1}). From (\ref{b2-rel-1}%
) and (\ref{D,C+D-commutativity}) it follows that $\Delta _{\widehat{ch}%
}^{R}\left( C_{t}+D_{t}\right) \mid _{t\neq s}$ are also invariant:
\begin{eqnarray}
\mathcal{F}_{s}^{\widehat{J}}\circ \left( \Delta _{\widehat{ch}}^{R}\left(
C_{t}+D_{t}\right) \right)  &\mid &_{s\neq t}=  \notag \\
&=&e^{\mathrm{ad}H_{s}^{\perp }\otimes \ln \left( C_{s}+D_{s}\right) }\circ
\left( D_{t}\otimes 1+C_{t}\otimes \left( C_{t}+D_{t}\right) \right) =
\notag \\
&=&D_{t}\otimes 1+C_{t}\otimes \left( C_{t}+D_{t}\right) .
\end{eqnarray}
Thus the Hopf subalgebra $\frak{B}_{t}\mid _{t\neq s}$ is invariant with
respect to the twisting by $\mathcal{F}_{s}^{\widehat{J}}$. The latter means
that the product of all the factors $\mathcal{F}_{s}^{\widehat{J}}$
\begin{equation}
\mathcal{F}^{\widehat{J}}=\prod_{s=1}^{p}\mathcal{F}_{s}^{\widehat{J}%
}=\prod_{s=1}^{p}e^{H_{s}^{\widehat{\perp }}\otimes \ln \left(
C_{s}+D_{s}\right) }=\prod_{s=1}^{p}e^{H_{s}^{\widehat{\perp }}\otimes
\omega _{s}}
\end{equation}
can be applied to $U_{\widehat{ch}}^{R}\left( A_{n-1}\right) $
\begin{equation}
\mathcal{F}^{\widehat{J}}:U_{\widehat{ch}}^{R}\left( A_{n-1}\right)
\longrightarrow U_{\mathcal{P}}\left( A_{n-1}\right) .
\end{equation}
The deformed algebra $U_{\mathcal{P}}\left( A_{n-1}\right) $ can be also
presented as twisted by the product of the rotated chain and the
quasi-Jordanian factors
\begin{equation}
\mathcal{F}_{\mathcal{P}}=\mathcal{F}^{\widehat{J}}\mathcal{F}_{\widehat{ch}%
}^{R},
\end{equation}
\begin{equation}
\mathcal{F}_{\mathcal{P}}:U\left( A_{n-1}\right) \longrightarrow U_{\mathcal{%
P}}\left( A_{n-1}\right) .
\end{equation}
The carrier of the twist $\mathcal{F}_{\mathcal{P}}$,

\begin{equation}
g_{\mathcal{P}}=\left( \bigoplus_{s=1}^{p}sl\left( 2\right) _{i}\right)
\vartriangleright \left( g_{ch}\setminus \sum_{s=1}^{p}g^{\gamma
_{s}}\right), 
\end{equation}
is the parabolic subalgebra.

\section{Examples.}

\subsection{$sl\left( 4\right) $}

Here the chain contains two links
\begin{equation*}
\mathcal{F}_{\widehat{ch}}=e^{\widehat{H_{2}}\otimes \sigma
_{2,3}}e^{E_{1,2}\otimes E_{2,4}e^{-\sigma _{1,4}}}e^{E_{1,3}\otimes
E_{3,4}e^{-\sigma _{1,4}}}\cdot e^{\widehat{H_{1}}\otimes \sigma _{1,4}}
\end{equation*}
with
\begin{eqnarray}
\widehat{H_{1}} &=&\frac{1}{4}\mathbf{I}-E_{4,4}; \\
\widehat{H_{2}} &=&-\frac{3}{4}\mathbf{I}+E_{1,1}+E_{2,2}+E_{4,4}.
\end{eqnarray}

The main parameters of this chain deformation are

\begin{eqnarray}
m &=&2;\quad r=3; \\
p &=&\dim H^{\perp }=1;
\end{eqnarray}
The generators included in the chain carrier subalgebra and the possible
choice of external coordinates $\mathbf{E}_{k+1,k}$,
\begin{equation}
\left[
\begin{array}{llll}
E_{1,1} & E_{1,2} & E_{1,3} & \fbox{E$_{1,4}$} \\
\mathbf{E}_{2,1} & E_{2,2} & \fbox{E$_{2,3}$} & E_{2,4} \\
&  & E_{3,3} & E_{3,4} \\
&  &  & E_{4,4}
\end{array}
\right] \quad \mathrm{or}\quad \left[
\begin{array}{llll}
E_{1,1} & E_{1,2} & E_{1,3} & \fbox{E$_{1,4}$} \\
& E_{2,2} & \fbox{E$_{2,3}$} & E_{2,4} \\
&  & E_{3,3} & E_{3,4} \\
&  & \mathbf{E}_{4,3} & E_{4,4}
\end{array}
\right] ,
\end{equation}
show us that there are two minimal parabolic subalgebras of the type $g_{%
\mathcal{P}}$ appropriate for our twist algorithm. They are defined by the
sets
\begin{equation}
\begin{array}{l}
\Psi _{\mathcal{P}_{1}}=\left\{ \alpha _{2},\alpha _{3},\right\} , \\
\Psi _{\mathcal{P}_{3}}=\left\{ \alpha _{1},\alpha _{2},\right\} ,
\end{array}
\qquad \alpha _{k}:=e_{k}-e_{k+1},
\end{equation}
\begin{equation}
\begin{array}{l}
\Gamma _{\mathcal{P}_{1}}=\left\{ \alpha _{j\left( k\right) }\mid
k=1,j\left( 1\right) =4-\chi \left( 1\right) =1\right\} =\left\{ \alpha
_{1}\right\} , \\
\Gamma _{\mathcal{P}_{3}}=\left\{ \alpha _{j\left( k\right) }\mid
k=3,j\left( 3\right) =4-\chi \left( 3\right) =3\right\} =\left\{ \alpha
_{3}\right\} ,
\end{array}
\end{equation}
\begin{equation}
\begin{array}{l}
\Phi \left( \Psi _{\mathcal{P}_{1}}\right) =\Lambda ^{+}\bigcup \left\{
-\alpha _{1}\right\} , \\
\Phi \left( \Psi _{\mathcal{P}_{3}}\right) =\Lambda ^{+}\bigcup \left\{
-\alpha _{3}\right\} .
\end{array}
\end{equation}
and can be presented in the following form:
\begin{equation}
g_{\mathcal{P}_{1}}=H+\sum_{\lambda \in \Lambda ^{+}}g^{\lambda }+g^{-\alpha
_{1}}=B^{+}+g^{-\alpha _{1}},
\end{equation}
\begin{equation}
g_{\mathcal{P}_{3}}=H+\sum_{\lambda \in \Lambda ^{+}}g^{\lambda }+g^{-\alpha
_{3}}=B^{+}+g^{-\alpha _{3}}.
\end{equation}

First let us study the case $g_{\mathcal{P}_{1}}$. In the one-dimensional
space $H^{\perp }$ take the generator

\begin{equation}
H_{1}^{\widehat{\perp }}=\frac{1}{2}\mathbf{I}-\left( E_{1,1}+E_{4,4}\right)
.
\end{equation}
and consider the external coordinate:

\begin{equation*}
\widehat{E}_{1}=E_{2,1}+E_{3,4}e^{\sigma _{2,3}-\sigma _{1,4}};
\end{equation*}

The peculiarity of the Cartan subalgebras in $sl\left( 4k\right) $ , $k\in
\mathsf{N}^{+}$ , is that the following relation holds
\begin{equation}
2H_{n-3,n-2}+\widehat{H_{n-2}}-\widehat{H_{n-3}}=0.  \label{4k-property}
\end{equation}
Thus in the twisted coproduct $\Delta _{\widehat{ch}}\left( \widehat{E}%
_{1}\right) $ the term connecting the space $H^{\perp }$ and $\widehat{E}_{1}
$ is zero and both the Cartan generator $H_{1}^{\perp }$ and the external
coordinate have primitive coproducts:
\begin{eqnarray}
D_{1} &=&\widehat{E}_{1}, \\
\Delta _{\widehat{ch}}\left( D_{1}\right)  &=&D_{1}\otimes 1+1\otimes D_{1};
\\
\Delta _{\widehat{ch}}\left( H_{1}^{\perp }\right)  &=&H_{1}^{\perp }\otimes
1+1\otimes H_{1}^{\perp }.
\end{eqnarray}
It is easy to check that they obey the $b^{2}$-relations
\begin{equation}
\left[ H_{1}^{\perp },D_{1}\right] =D_{1}.
\end{equation}
These are the standard conditions guaranteeing that the ordinary Jordanian
twist
\begin{equation}
\mathcal{F}_{_{1}}^{\mathcal{J}}=e^{H_{1}^{\perp }\otimes \ln \left(
1+D_{1}\right) }=e^{H_{1}^{\perp }\otimes \ln \left( 1+\widehat{E}%
_{1}\right) }=e^{H_{1}^{\perp }\otimes \omega _{1}}
\end{equation}
is a solution of the twist equations for $\Delta _{\widehat{ch}}$. As a
consequence the product
\begin{equation}
\mathcal{F}_{\mathcal{P}_{1}}=\mathcal{F}_{_{1}}^{\mathcal{J}}\mathcal{F}_{%
\widehat{ch}}
\end{equation}
is a solution of the twist equations for the undeformed $U\left( sl\left(
4\right) \right) ,$
\begin{equation}
\mathcal{F}_{\mathcal{P}_{1}}:U\left( sl\left( 4\right) \right)
\longrightarrow U_{\mathcal{P}_{1}}\left( sl\left( 4\right) \right) ,
\end{equation}
with the parabolic carrier $g_{\mathcal{P}_{1}}$.

The alternative possibility is to choose $g_{\mathcal{P}_{3}}$ as the
carrier subalgebra. In this case taking into account (\ref{4k-property}) the
external coordinate can be presented as
\begin{eqnarray}
\widehat{E}_{3} &=&E_{4,3}-\left( \widehat{H_{1}}-\widehat{H_{2}}\right)
E_{1,3}+E_{1,2}=  \notag \\
&=&E_{4,3}-2H_{3,4}E_{1,3}+E_{1,2}.
\end{eqnarray}
The direct computation shows that the twisted coproduct $\Delta _{\widehat{ch%
}}\left( E_{3}^{\symbol{94}}\right) $\ again has no terms like $H^{\perp
}\bigwedge E$ . This time it is quasiprimitive:
\begin{eqnarray}
D_{3} &=&\widehat{E}_{3}e^{\sigma _{1,4}-\sigma _{2,3}},\quad
C_{3}=e^{\sigma _{1,4}-\sigma _{2,3}} \\
\Delta _{\widehat{ch}}\left( D_{3}\right)  &=&D_{3}\otimes 1+C_{3}\otimes
D_{3}; \\
\Delta _{\widehat{ch}}\left( H_{1}^{\perp }\right)  &=&H_{1}^{\perp }\otimes
1+1\otimes H_{1}^{\perp }.
\end{eqnarray}
The pair $\left\{ H_{1}^{\perp },D_{3}\right\} $ forms a Borel subalgebra:
\begin{equation}
\left[ H_{1}^{\perp },D_{3}\right] =-D_{3};
\end{equation}
and according to Lemma 4 the factor
\begin{equation}
\mathcal{F}_{_{3}}^{\widehat{J}}=e^{-H_{1}^{\perp }\otimes \ln \left(
C_{3}+D_{3}\right) }=e^{-H_{1}^{\perp }\otimes \ln \left( \left( 1+\widehat{E%
}_{3}\right) e^{\sigma _{1,4}-\sigma _{2,3}}\right) }=e^{-H_{1}^{\perp
}\otimes \omega _{3}}
\end{equation}
is a solution to the twist equations for $\Delta _{\widehat{ch}}.$ The same
is true for the product
\begin{equation}
\mathcal{F}_{\mathcal{P}_{3}}=\mathcal{F}_{_{3}}^{\widehat{J}}\mathcal{F}_{%
\widehat{ch}}
\end{equation}
with respect to the undeformed $\Delta $ and thus $\mathcal{F}_{\mathcal{P}%
_{3}}$ represents the parabolic twist with the carrier $g_{\mathcal{P}_{3}}$.

Notice that these twists do not commute and we cannot apply both
deformations simultaneously, no combinations of $\omega _{1}$ and $\omega
_{3}$ are allowed.

\subsection{$sl\left( 11\right) $}

For this case we have
\begin{eqnarray}
m &=&5;\quad r=10; \\
p &=&\dim H^{\perp }=5;
\end{eqnarray}
\begin{eqnarray}
\Psi _{\mathcal{P}} &=&\left\{ \alpha _{\chi \left( s\right) }\mid \chi
\left( s\right) =\frac{1}{2}\left( 11+\left( 11-2s\right) \left( -1\right)
^{s+1}\right) ,s=1,\ldots ,5\right\} =  \notag \\
&=&\left\{ \alpha _{2},\alpha _{4},\alpha _{6},\alpha _{8},\alpha
_{10}\right\} ,\qquad \alpha _{k}:=e_{k}-e_{k+1},
\end{eqnarray}
\begin{equation}
\Phi \left( \Psi \right) =\Lambda ^{+}\bigcup \left\{ -\alpha _{j}\right\}
,\qquad \alpha _{j}\in \Gamma _{\mathcal{P}},
\end{equation}
where
\begin{equation}
\Gamma _{\mathcal{P}}=\left\{ \alpha _{j}\mid j\left( s\right) =11-\chi
\left( s\right) \right\} =\left\{ \alpha _{1},\alpha _{3},\alpha _{5},\alpha
_{7},\alpha _{9}\right\} .
\end{equation}
Consider the parabolic subalgebra
\begin{equation}
g_{\mathcal{P}}=H+\sum_{\lambda \in \Lambda ^{+}}g^{\lambda }+\sum_{\gamma
\in \Gamma _{\mathcal{P}}}g^{-\gamma },
\end{equation}
with the set of 2-dimensional Borel subalgebras
\begin{equation}
b_{s}=\left\{ H_{s}^{\widehat{\perp }},g^{-\alpha _{j\left( s\right)
}}\right\} _{s=1,\ldots ,5},
\end{equation}
where
\begin{eqnarray}
H_{1}^{\widehat{\perp }} &=&\frac{2}{11}\mathbf{I}-\left(
E_{1,1}+E_{11,11}\right) ;  \notag \\
H_{2}^{\widehat{\perp }} &=&-\frac{4}{11}\mathbf{I}+\sum_{u=1}^{2}\left(
E_{u,u}+E_{12-u,12-u}\right) ;  \notag \\
H_{3}^{\widehat{\perp }} &=&\frac{6}{11}\mathbf{I}-\sum_{u=1}^{3}\left(
E_{u,u}+E_{12-u,12-u}\right) ;  \notag \\
H_{4}^{\widehat{\perp }} &=&-\frac{8}{11}\mathbf{I}+\sum_{u=1}^{4}\left(
E_{u,u}+E_{12-u,12-u}\right) ;  \notag \\
H_{5}^{\widehat{\perp }} &=&\frac{10}{11}\mathbf{I}-\sum_{u=1}^{5}\left(
E_{u,u}+E_{12-u,12-u}\right) .
\end{eqnarray}

The dual coordinates $\widehat{E}_{s}$ for the generators of the root
subspaces $g^{-\alpha _{j\left( s\right) }}$ are as follows:
\begin{eqnarray}
\widehat{E}_{1} &=&E_{2,1}+E_{10,11}e^{\sigma _{2,10}-\sigma _{1,11}};
\notag \\
\widehat{E}_{2} &=&E_{10,9}+\left( \widehat{H}_{3}-\widehat{H}_{2}\right)
E_{2,9}+E_{2,3};  \notag \\
\widehat{E}_{3} &=&E_{4,3}+E_{8,9}e^{\sigma _{4,8}-\sigma _{3,9}};  \notag \\
\widehat{E}_{4} &=&E_{8,7}+\left( \widehat{H}_{5}-\widehat{H}_{4}\right)
E_{4,7}+E_{4,5};  \notag \\
\widehat{E}_{5} &=&E_{6,5};
\end{eqnarray}

After the application of the full chain twist
\begin{eqnarray*}
\mathcal{F}_{\widehat{ch}} &=&\prod_{l=1}^{5}\mathcal{F}_{\widehat{\mathrm{%
link}}}^{l} \\
&=&e^{E_{5,6}\otimes E_{6,7}e^{-\sigma _{5,7}}}e^{\widehat{H_{5}}\otimes
\sigma _{5,7}}\cdot  \\
&&\cdot e^{\sum_{k_{4}=5}^{7}E_{4,k_{4}}\otimes E_{k_{4},8}}e^{\widehat{H_{4}%
}\otimes \sigma _{4,8}}\cdot  \\
&&\cdot e^{\sum_{k_{3}=4}^{8}E_{3,k_{3}}\otimes E_{k_{3},9}e^{-\sigma
_{3,9}}}e^{\widehat{H_{3}}\otimes \sigma _{3,9}}\cdot  \\
&&\cdot e^{\sum_{k_{2}=3}^{9}E_{2,k_{2}}\otimes E_{k_{2},10}}e^{\widehat{%
H_{2}}\otimes \sigma _{2,10}}\cdot  \\
&&\cdot e^{\sum_{k_{1}=2}^{10}E_{1,k_{1}}\otimes E_{k_{1},11}e^{-\sigma
_{1,11}}}e^{\widehat{H_{1}}\otimes \sigma _{1,11}}
\end{eqnarray*}

with
\begin{equation}
\begin{array}{l}
\widehat{H}_{1}=+\frac{1}{11}\mathbf{I}-E_{11,11}; \\
\widehat{H}_{2}=-\frac{3}{11}\mathbf{I}+E_{1,1}+E_{2,2}+E_{11,11}; \\
\widehat{H}_{3}=+\frac{5}{11}\mathbf{I}%
-E_{1,1}-E_{2,2}-E_{9,9}-E_{10,10}-E_{11,11}; \\
\widehat{H}_{4}=-\frac{7}{11}\mathbf{I}+E_{1,1}+\ldots
+E_{4,4}+E_{9,9}+E_{10,10}+E_{11,11}; \\
\widehat{H}_{5}=+\frac{9}{11}\mathbf{I}-E_{1,1}-\ldots
-E_{4,4}-E_{7,7}-E_{8,8}-\ldots -E_{11,11};
\end{array}
\end{equation}
to $U\left( sl\left( 11\right) \right) $ we obtain the following coproduct
maps for the external coordinates
\begin{eqnarray}
\Delta _{\widehat{ch}}\left( \widehat{E}_{1}\right)  &=&\widehat{E}%
_{1}\otimes 1+1\otimes \widehat{E}_{1}+\left( -2H_{1}^{\widehat{\perp }%
}\right) \otimes E_{2,11}e^{-\sigma _{1,11}};  \notag \\
\Delta _{\widehat{ch}}\left( \widehat{E}_{2}\right)  &=&\widehat{E}%
_{2}\otimes 1+1\otimes \widehat{E}_{2}+E_{2,9}\otimes \left( -2H_{2}^{%
\widehat{\perp }}\right) ;  \notag \\
\Delta _{\widehat{ch}}\left( \widehat{E}_{3}\right)  &=&\widehat{E}%
_{3}\otimes 1+1\otimes \widehat{E}_{3}+\left( -2H_{3}^{\widehat{\perp }%
}\right) \otimes E_{4,9}e^{-\sigma _{3,9}};  \notag \\
\Delta _{\widehat{ch}}\left( \widehat{E}_{4}\right)  &=&\widehat{E}%
_{4}\otimes 1+1\otimes \widehat{E}_{4}+E_{4,7}\otimes \left( -2H_{4}^{%
\widehat{\perp }}\right) ;  \notag \\
\Delta _{\widehat{ch}}\left( \widehat{E}_{5}\right)  &=&\widehat{E}%
_{5}\otimes 1+1\otimes \widehat{E}_{5}+\left( -2H_{5}^{\widehat{\perp }%
}\right) \otimes E_{6,7}e^{-\sigma _{5,7}};
\end{eqnarray}
while the Cartan generators $H_{s}^{\widehat{\perp }}$ remain primitive:
\begin{equation}
\Delta _{\widehat{ch}}\left( H_{s}^{\widehat{\perp }}\right) =H_{s}^{%
\widehat{\perp }}\otimes 1+1\otimes H_{s}^{\widehat{\perp }}.
\end{equation}
Let us introduce the set
\begin{equation}
C_{s}=e^{2\sum_{i=1}^{s}\sigma _{i,n-i+1}};
\end{equation}
Taking into account the commutation properties of $\sigma $\ 's with the
external coordinates $\widehat{E}_{s}$:

\begin{equation*}
\begin{array}{llllll}
& \widehat{E}_{1} & \widehat{E}_{2} & \widehat{E}_{3} & \widehat{E}_{4} &
\widehat{E}_{5} \\
\sigma _{1,11} & -E_{2,11}e^{-\sigma _{1,11}} & 0 & 0 & 0 & 0 \\
\sigma _{2,10} & +E_{2,11}e^{-\sigma _{1,11}} & +E_{2,9} & 0 & 0 & 0 \\
\sigma _{3,9} & 0 & -E_{2,9} & -E_{4,9}e^{-\sigma _{3,9}} & 0 & 0 \\
\sigma _{4,8} & 0 & 0 & +E_{4,9}e^{-\sigma _{3,9}} & +E_{4,7} & 0 \\
\sigma _{5,7} & 0 & 0 & 0 & -E_{4,7} & -E_{6,7}e^{-\sigma _{5,7}}
\end{array}
\end{equation*}
we obtain the main set of commutation relations:

\begin{equation*}
\left( -1\right) ^{s}\left[ \ln \left( C_{s}\right) ,\widehat{E}_{s}\right]
=2\delta _{s,t}\cdot \left\{
\begin{array}{c}
E_{s+1,n-s+1}e^{-\sigma _{s,n-s+1}};\quad \text{for odd }s \\
E_{s,n-s};\quad \text{for even }s
\end{array}
\right.
\end{equation*}
\begin{eqnarray*}
\left[ \ln \left( C_{t}\right) ,E_{s+1,n-s+1}e^{-\sigma _{s,n-s+1}}\right]
&=&0; \\
\left[ \ln \left( C_{t}\right) ,E_{s,n-s}\right]  &=&0; \\
\left[ C_{s},C_{t}\right]  &=&0.
\end{eqnarray*}
Compose the rotation twist:
\begin{equation}
\mathcal{F}^{R}=\exp \left(
\begin{array}{c}
-H_{1}^{\widehat{\perp }}\otimes \ln \left( C_{1}\right) +\ln \left(
C_{2}\right) \otimes H_{2}^{\widehat{\perp }}-H_{3}^{\widehat{\perp }%
}\otimes \ln \left( C_{3}\right)  \\
+\ln \left( C_{4}\right) \otimes H_{4}^{\widehat{\perp }}-H_{5}^{\widehat{%
\perp }}\otimes \ln \left( C_{5}\right)
\end{array}
\right)
\end{equation}
and apply it to the algebra $U_{\widehat{ch}}\left( sl\left( 11\right)
\right) $
\begin{equation}
\mathcal{F}^{R}:U_{\widehat{ch}}\left( sl\left( 11\right)\right) 
\longrightarrow U_{%
\widehat{ch}}^{R}\left( sl\left( 11\right)\right)
\end{equation}
the coproducts of the external coordinates become quasiprimitive:
\begin{equation}
\Delta _{\widehat{ch}}^{R}\left( \widehat{E}_{s}\right) =\left\{
\begin{array}{c}
\widehat{E}_{s}\otimes C_{s}^{-1}+1\otimes \widehat{E}_{s},\text{\quad }%
s=2k-1, \\
\widehat{E}_{s}\otimes 1+C_{s}\otimes \widehat{E}_{s},\quad s=2k.
\end{array}
\right. ;
\end{equation}
while the coproducts $\Delta _{\widehat{ch}}^{R}\left( H_{s}^{\widehat{\perp
}}\right) $ remain primitive. In terms of redefined external coordinates
\begin{equation}
D_{s}=\left\{
\begin{array}{c}
\left( E_{2,1}+E_{10,11}e^{\sigma _{2,10}-\sigma _{1,11}}\right) e^{2\sigma
_{1,11}}, \\
E_{10,9}+\left( \widehat{H}_{3}-\widehat{H}_{2}\right) E_{2,9}+E_{2,3}, \\
\left( E_{4,3}+E_{8,9}e^{\sigma _{4,8}-\sigma _{3,9}}\right)
e^{2\sum_{i=1}^{3}\sigma _{i,n-i+1}}, \\
E_{8,7}+\left( \widehat{H}_{5}-\widehat{H}_{4}\right) E_{4,7}+E_{4,5}, \\
E_{6,5}e^{2\sum_{i=1}^{5}\sigma _{i,n-i+1}}.
\end{array}
\right.
\end{equation}
these coproducts are uninform:
\begin{equation}
\Delta _{\widehat{ch}}^{R}\left( D_{s}\right) =D_{s}\otimes 1+C_{s}\otimes
D_{s}.
\end{equation}
\qquad

For any $s=1,\ldots ,p$ the subalgebras $\frak{B}_{s}$ generated by $\left\{
H_{s}^{\widehat{\perp }},D_{s},C_{s}\right\} $ comply the conditions of
Lemma 4 and the twisting elements
\begin{equation}
\mathcal{F}_{s}^{\widehat{J}}=e^{H_{s}^{\widehat{\perp }}\otimes \ln \left(
C_{s}+D_{s}\right) }=e^{H_{s}^{\widehat{\perp }}\otimes \omega _{s}}
\end{equation}
with
\begin{equation*}
\omega _{s}=\ln \left( C_{s}+D_{s}\right)
\end{equation*}
are the solutions to the Drinfeld equations (\ref{twist-equations}). The
product
\begin{equation}
\mathcal{F}^{\widehat{J}}=\prod_{s=1}^{5}\mathcal{F}_{s}^{\widehat{J}%
}=\prod_{s=1}^{5}e^{H_{s}^{\widehat{\perp }}\otimes \ln \left(
C_{s}+D_{s}\right) }=\prod_{s=1}^{5}e^{H_{s}^{\widehat{\perp }}\otimes
\omega _{s}};
\end{equation}
can be applied to $U_{\widehat{ch}}^{R}\left( sl\left( 11\right) \right) $
\begin{equation}
\mathcal{F}^{\widehat{J}}:U_{\widehat{ch}}^{R}\left( sl\left( 11\right)
\right) \longrightarrow U_{\mathcal{P}}\left( sl\left( 11\right) \right) .
\end{equation}
The product
\begin{equation}
\mathcal{F}_{\mathcal{P}}=\mathcal{F}^{\widehat{J}}\mathcal{F}_{\widehat{ch}%
}^{R},
\end{equation}
\begin{equation}
\mathcal{F}_{\mathcal{P}}:U\left( sl\left( 11\right) \right) \longrightarrow
U_{\mathcal{P}}\left( sl\left( 11\right) \right)
\end{equation}
has the parabolic carrier

\begin{equation}
g_{\mathcal{P}}=\left( \bigoplus_{s=1}^{5}sl\left( 2\right) _{s}\right)
\vartriangleright \left( g_{ch}\setminus \sum_{s=1}^{5}g^{\gamma
_{s}}\right) .
\end{equation}

\section{Conclusions}

We have demonstrated that for any algebra $A_{n-1}$ ($n>2$) there exists a
twist deformation $\mathcal{F}_{\mathcal{P}}:U\left( A_{n-1}\right)
\longrightarrow U_{\mathcal{P}}\left( A_{n-1}\right) $ whose carrier $g_{%
\mathcal{P}}$ is a parabolic subalgebra with non-Abelian Levi factor. These
carriers
are the non-Abelian extensions of $g_{\widehat{ch}}$ that contain $\mathbf{B}%
^{+}(g)$ and have nontrivial intersection $g_{\mathcal{P}}\cap \mathbf{N}%
_{g}^{-}=\emptyset $ .

The corresponding twisting elements $\mathcal{F}_{\mathcal{P}}$ are the
products of deformed Jordanian twists $\mathcal{F}^{\widehat{J}}$ and the
rotated full chains $\mathcal{F}_{\widehat{ch}}^{R}$. According to the
quantum duality the Hopf subalgebras $U_{\mathcal{P}}\left( g_{\mathcal{P}%
}\right) $ are the quantized algebras $\mathrm{Fun}_{\mathcal{P}}\left( G_{%
\mathcal{P}}\right) $ of coordinate functions on the universal covering
groups $G_{\mathcal{P}}$ (with Lie algebras $g_{\mathcal{P}}$). Notice that
in this type of quantization the deformation of algebra $\mathrm{Fun}\left(
G_{\mathcal{P}}\right) $ is due to the fact that the coordinates are subject
to the commutation relations of $g_{\mathcal{P}}$. Thus we have constructed
quantum groups for a class of parabolic subgroups of linear groups $SL\left(
n\right) $.

It was mentioned in the Introduction that for series $B_{n}$,$C_{n}$ and $%
D_{n=2k}$ the full chains have (minimal) parabolic carriers. In the case of
even-odd orthogonal algebras $D_{n=2k+1}$ the chain carrier is not
parabolic, $\mathrm{dim}\left( \mathbf{B}^{+}(g)\setminus g_{ch}\right) >0$.
Due to the canonical isomorphism one of such algebras has been already
treated in the section Examples where we considered $g=sl\left( 4\right)
\approx so\left( 6\right) $. Thus it is clear that the problem of
constructing parabolic twists for $D_{2k+1}$ can be solved using the same
approach as in the above study.

It can be shown that the parabolic twists $\mathcal{F}_{\mathcal{P}}=%
\mathcal{F}^{\widehat{J}}\mathcal{F}_{\widehat{ch}}^{R}$ are completely
factorizable. The factors in the decomposition $\mathcal{F}_{\mathcal{P}%
}=\prod_{s=1}^{p}\mathcal{F}_{s}^{\widehat{J}}\mathcal{F}^{R}\prod_{l=1}^{m}
\mathcal{F}_{\widehat{link}}^{l}$ can be supplied 
with the independent deformation
parameters $\left\{ \chi _{i},\psi _{j}\mid i=1,\ldots ,m;\right.$
$ \left. j=1,\ldots p\right\} $ and any subset of the factors $\mathcal{F}_{s}^{%
\widehat{J}}$ and $\mathcal{F}_{\widehat{lk}}^{lR}$ can be switched off by
tending the corresponding set of parameters to their limit values. The
factors $\mathcal{F}_{s}^{\widehat{J}}$ depend on the parameters of the
chain $\mathcal{F}_{\widehat{ch}}^{R}$. It is essential to mention that when
the whole chain is switched off the factors $\mathcal{F}_{s}^{\widehat{J}}$
acquire the canonical Jordanian form. The important consequence of this
property is that the parabolic twists $\mathcal{F}_{\mathcal{P}}$ are
products of the same set of basic twisting factors\textbf{\ }\cite{L-BTF} as
used in different forms of chains. Thus parabolic twists can also be
considered as chains of the special type.

When the variables $\left\{ \chi _{i},\psi _{j}\right\} $ are proportional
to the overall deformation parameter $\tau $ ($\chi _{i}=\tau \xi _{i},\psi
_{j}=\tau \zeta _{j}$) the first order terms in the expansion of the $%
\mathcal{R}$-matrices $\mathcal{R}_{\mathcal{P}}\left( \tau \right) =%
\mathcal{F}_{\mathcal{P}21}\left( \tau \right) \mathcal{F}_{\mathcal{P}%
}^{-1}\left( \tau \right) $\ give the parabolic classical $r$-matrices $%
r_{ech}$ constructed in \cite{LL}. Twisting by $\mathcal{F}_{\mathcal{P}}$
performs the quantization of $r_{ech}$.

\section{Acknowledgements}

The author is grateful to Prof. P.P.Kulish for stimulating discussions. The
work was supported by the Russian Foundation for Fundamental Research, grant
N 030100593.

\end{document}